\documentclass[11pt,reqno,tbtags,a4paper]{amsart}
\usepackage{amssymb}
\usepackage{url}
\usepackage[square,numbers]{natbib}
\bibpunct[, ]{[}{]}{;}{n}{,}{,}

\title
{A graphon counter example}

\date{6 September, 2019}

\author{Svante Janson}
\thanks{Partly supported by the Knut and Alice Wallenberg Foundation}
\address{Department of Mathematics, Uppsala University, PO Box 480,
SE-751~06 Uppsala, Sweden}
\email{svante.janson@math.uu.se}
\urladdr{http://www.math.uu.se/svante-janson}

\subjclass[2010]{} 

\numberwithin{equation}{section}

\renewcommand\le{\leqslant}
\renewcommand\ge{\geqslant}

\allowdisplaybreaks

\theoremstyle{plain}
\newtheorem{theorem}{Theorem}

\newtheorem*{problem}{Problem}

\theoremstyle{definition}

\newtheorem{remark}[theorem]{Remark}

\theoremstyle{remark}

\newenvironment{romenumerate}[1][-10pt]{
\addtolength{\leftmargini}{#1}\begin{enumerate}
 \renewcommand{\labelenumi}{\textup{(\roman{enumi})}}%
 }{\end{enumerate}}

\newcounter{oldenumi}
{\setcounter{oldenumi}{\value{enumi}}
\begin{romenumerate} \setcounter{enumi}{\value{oldenumi}}}
{\end{romenumerate}}

\newcounter{thmenumerate}

\newcounter{xenumerate}   

\newcommand\marginal[1]{\marginpar[\raggedleft\tiny #1]{\raggedright\tiny#1}}
\newcommand\REM[1]{{\raggedright\texttt{[#1]}\par\marginal{XXX}}}
\newcommand\XREM[1]{\relax}

\begingroup
  \divide\count255 by 60
  \multiply\count255 by -60
  \advance\count255 by \time
  \ifnum \count255 < 10 \xdef\klockan{\the\count1.0\the\count255}
  \else\xdef\klockan{\the\count1.\the\count255}\fi
\endgroup

\newcommand\set[1]{\ensuremath{\{#1\}}}
\newcommand\bigset[1]{\ensuremath{\bigl\{#1\bigr\}}}

\newcommand\xpar[1]{(#1)}
\newcommand\bigpar[1]{\bigl(#1\bigr)}
\newcommand\Bigpar[1]{\Bigl(#1\Bigr)}

\newcommand\Bigcpar[1]{\Bigl\{#1\Bigr\}}

\def\rompar(#1){\textup(#1\textup)}    

\def\xexp(#1){e^{#1}}

\newcommand\punkt{.\spacefactor=1000}    
    
\newcommand\ie{i.e\punkt}
\newcommand\eg{e.g\punkt}

\newcommand{\aex}{a.e\punkt}

\newcounter{CC}
\newcounter{cc}

\newcommand\gf{\varphi}

\newcommand\gl{\lambda}

\newcommand\gO{\Omega}
\newcommand\gs{\sigma}

\newcommand\eps{\varepsilon}

\renewcommand\phi{\xxx}  

\newcommand\cF{\mathcal F}

\newcommand\Bigindic[1]{\boldsymbol1\Bigcpar{#1}} 

\newcommand\qw{^{-1}}

\newcommand\intoi{\int_0^1}

\newcommand\oi{\ensuremath{[0,1]}}

\newcommand\setoi{\set{0,1}}

\newcommand\dd{\,\mathrm{d}}

\newcommand\leb{\gl}
\newcommand\mpr{measure preserving}
\newcommand\DD{\mathfrak D} 

\begin{document}

\begin{abstract} 
We give an example of a graphon such that there is no equivalent graphon
with a degree function that is (weakly) increasing.
\end{abstract}

\maketitle

\section{Introduction}\label{S:intro}

A central fact in the theory of graph limits 
(see \eg{} the book by \citet{Lovasz})
is that each graph limit can be represented by a graphon, but this
representation is not unique.
We say that two graphons are \emph{equivalent} if they define the same graph
limit; 
thus there is a
bijection between graph limits and equivalence classes of graphons.

Recall that graphons are symmetric measurable functions
$W:\gO\times\gO\to\oi$, where $\gO=(\gO,\cF,\mu)$ is a probability space.
We may always choose $\gO$ to be $\oi$ with Lebesgue measure,
in the sense that any graphon is equivalent to a graphon defined on $\oi$,
but it is often advantageous to use graphons defined on other 
probability spaces $\gO$ too.

The characterization of equivalence between graphons is known to be 
complicated; it includes \aex{} equality and taking the pull-back 
by a \mpr{} map
(see below for definitions), 
but is not limited to this. See \eg{}
\cite{LovaszSz}, \cite{BR}, \cite{SJ209}, \cite{BCL} and \cite{SJ249}.

Given a graph limit, it would be desirable to somehow define a canonical
graphon representing it
(at least up to \aex{} equality); 
in other words, to define a canonical choice of a graphon in the
corresponding equivalence class.
In some special cases, this can be done in a natural way. 
For example, 
see \cite{SJ238},
a graph limit that is the limit of a sequence of threshold
graphs
can always be represented by a graphon $W(x,y)$
on $\oi$ that only takes values in
\setoi, and furthermore is increasing 
in each coordinate separately
(we say that a function 
$f(x)$ is increasing if $f(x)\le f(y)$ when $x\le y$); moreover,
two such graphons are equivalent if and only if they are \aex{} equal.
There is thus a canonical graphon representing each threshold
graph limit.

Similarly, if a graphon $W(x,y)$ defined on $\oi$ has  a degree function
\begin{equation}\label{degree}
\DD(x)=  \DD _W(x):=\intoi W(x,y)\dd y
\end{equation}
that is a strictly increasing function $\oi\to\oi$,
then it not difficult to show 
that any equivalent graphon that also has an increasing degree function
is \aex{} equal to $W$.
(Use \eqref{bb} below. We  omit the details.)
Hence, a graphon with a strictly increasing degree function can be regarded
as a canonical choice in its equivalence class.

Of course, not every graphon is equivalent to such a graphon; for example
not a graphon with a constant degree function. Nevertheless, this leads to
the following interesting question.
We repeat that we use 'increasing' in the weak sense
(also known as 'weakly increasing'):
$f$ is increasing if $f(x)\le f(y)$ when $x\le y$; 

\begin{problem}
  Given any graphon $W$, does there exist an equivalent graphon on $\oi$
with an
increasing degree function \eqref{degree}?
\end{problem}

The purpose of this note is to show that this is \emph{not} the case.

\begin{theorem}
  There exists a graphon on $\oi$ such that there is no equivalent graphon
on $\oi$   with a  (weakly) increasing degree function.
\end{theorem}

We prove this theorem by giving a simple explicit  example in
\eqref{wessex}.
The example is similar to, and inspired by, standard examples such as
\cite[Example 7.11]{Lovasz}
showing that two equivalent graphons are not necessarily pull-backs of each
other.

\begin{remark}
The analogue for finite graphs of 
the problem above for graphons is the trivial fact that the vertices of a
graph can be ordered with (weakly) increasing vertex degrees.
Note that there will always be ties, so even for a finite graph, this does
not define a unique canonical labelling.
\end{remark}

\subsection{Some notation}

$\oi$ will, as above, be regarded as a probability space equipped with the
Lebesgue measure and the Lebesque $\gs$-field. 
(We might also use  the Borel $\gs$-field. 
For the present paper, this makes no difference; for other
purposes, the choice of $\gs$-field may have some importance.)

Let $(\gO_1,\cF_1,\mu_1)$  and $(\gO_2,\cF_2,\mu_2)$ be two probability
spaces. A function $\gf:\gO_1\to\gO_2$ is \emph{\mpr} if 
$\mu_1\xpar{\gf\qw(A)}=\mu_2(A)$ for any measurable $A\subseteq\gO_2$.
If $W$ is a graphon on $\gO_2$ and $\gf:\gO_1\to\gO_2$ is measure preserving, 
then the \emph{pull-back} $W^\gf$ is the graphon
$W^\gf(x,y):=W\bigpar{\gf(x),\gf(y)}$ defined on $\gO_1$.
As mentioned above, a pull-back $W^\gf$ is always equivalent to $W$.

\section{The example}

Our example is the graphon
\begin{align}\label{wessex}
  W(x,y):=
  \begin{cases}
    4xy, &x,y\in(0,\frac12), \\
    1/2, & x+y>3/2,\\
    0,& \text{otherwise}.
  \end{cases}
\end{align}
Note that   the degree function is given by
\begin{align}\label{d}
  \DD(x):=\intoi W(x,y)\dd y =
  \begin{cases}
    \frac12x, & x\in(0,\frac12),\\
    \frac12(x-\frac12),& x\in(\frac12,1).
  \end{cases}
\end{align}

Suppose that $W$ is equivalent to a graphon $W_1$ on $\oi$ that has an 
increasing degree function $\DD_1(x):=\intoi W_1(x,y)\dd y$;
we will show that this leads to a contradiction.

The equivalence $W\cong W_1$ implies by
\cite[Corollary 2.7]{BR}, see also \cite[Corollary 10.35]{Lovasz} and
  \cite[Theorem 8.6]{SJ249},
that there exist a probability space $(\gO,\mu)$ 
and two \mpr{} maps $\gf,\psi:\gO\to\oi$ such
that
$W^\gf=W_1^\psi$ a.e., \ie,
\begin{align}\label{a}
  W\bigpar{\gf(x),\gf(y)}
  =W_1\bigpar{\psi(x),\psi(y)},
  \qquad \text{a.e.\ on $\gO^2$}.
\end{align}
(The probability space $\gO,\mu)$ can be taken as
$\oi$ with Lebesgue measure, but we have no need for this.
Instead, we prefer to use the notation $\gO$ and $\mu$ to distinguish
between this space and $\oi$, which hopefully will 
make the proof easier to follow.)

Since $\gf$ and $\psi$ are \mpr, we have for every Borel measurable $f\ge0$
on $\oi$,
\begin{align}\label{b}
  \intoi f(x)\dd x = \int_\gO f(\gf(x))\dd\mu(x)
 = \int_\gO f(\psi(x))\dd\mu(x).
\end{align}
We use this repeatedly below.

In particular, \eqref{a} and \eqref{b} imply that for a.e.\ $x\in\gO$
\begin{align}\label{bb}
  \DD\bigpar{\gf(x)}
&  =\intoi W\bigpar{\gf(x),y}\dd y
  =\int_\gO W\bigpar{\gf(x),\gf(y)}\dd \mu(y)
\notag\\&  =\int_\gO W_1\bigpar{\psi(x),\psi(y)}\dd \mu(y)
  =\intoi W_1\bigpar{\psi(x),y}\dd y
  =\DD_1\bigpar{\psi(x)}.
\end{align}
Hence, for every real $r\in(0,\frac{1}4]$,
using \eqref{d},
\begin{align}
&  \leb\bigset{x\in\oi:\DD_1(x)\le r}
  =  \mu{\bigset{x\in\gO:\DD_1(\psi(x))\le r}}
  \notag\\&\qquad
=  \mu\bigset{x\in\gO:\DD(\gf(x))\le r}
  =
  \leb{\bigset{x\in\oi:\DD(x)\le r}}
  =4r.
\end{align}
Since we have assumed that $\DD_1$ is increasing, this implies
\begin{align}\label{d1}
  \DD_1(x)=x/4,
  \qquad x\in(0,1).
\end{align}

Define
\begin{align}
  \label{h}
  h(x):=\leb\bigset{y:W(x,y)\notin \set{0,\tfrac12}}
  =
  \begin{cases}
    \tfrac12, & x\in(0,\tfrac12),
    \\
    0,& x\in(\tfrac12,1),
  \end{cases}
\end{align}
and, similarly,
\begin{align}
  \label{h1}
  h_1(x):=\leb\bigset{y:W_1(x,y)\notin \set{0,\tfrac12}}.
\end{align}
Then \eqref{a} implies, similarly to \eqref{bb}, for \aex{} $x\in\gO$,
\begin{align}\label{hh}
  h(\gf(x))
 &=\leb\bigset{y:W(\gf(x),y)\notin\set{0,\tfrac12}}\notag\\
 &=\mu\bigset{y:W(\gf(x),\gf(y))\notin\set{0,\tfrac12}}\notag\\
&=    \mu\bigset{y:W_1(\psi(x),\psi(y))\notin\set{0,\tfrac12}}\notag\\
&=    \leb\bigset{y:W_1(\psi(x),y)\notin\set{0,\tfrac12}}
    =h_1\xpar{\psi(x)}.
\end{align}

This will yield our contradiction. We first calculate $h_1$.

If $0<a<b<1$, then, using 
\eqref{d1}, \eqref{b},
\eqref{hh},
\eqref{bb}, and \eqref{b} again,
\begin{align}
  \int_a^b h_1(x)\dd x
  &=\intoi h_1(x)\Bigindic{\frac{a}{4}<\DD_1(x)<\frac{b}{4}}\dd x
\notag\\&    
=\int_\gO h_1(\psi(x))\Bigindic{\frac{a}{4}<\DD_1(\psi(x))<\frac{b}{4}}\dd \mu(x)
\notag\\&    
  =\int_\gO h(\gf(x))\Bigindic{\frac{a}{4}<\DD(\gf(x))<\frac{b}{4}}\dd \mu(x)
\notag\\&    
=\intoi h(x)\Bigindic{\frac{a}{4}<\DD(x)<\frac{b}{4}}\dd x.\label{bx}
\end{align}
However, by \eqref{h} and \eqref{d},
\begin{align}
  \intoi h(x)\Bigindic{\frac{a}{4}<\DD(x)<\frac{b}{4}}\dd x&
=\frac12\int_0^{1/2}\Bigindic{\frac{a}{4}<\DD(x)<\frac{b}{4}}\dd x
                                                             \notag\\&    
=\frac12\leb\Bigpar{\frac{a}{2},\frac{b}{2}}=\frac{b-a}{4}.
\label{by}
\end{align}
Consequently, \eqref{bx} and \eqref{by} show that 
for every $a\in(0,1)$ and $\eps\in(0,1-a)$,
\begin{align}
  \frac{1}{\eps}\int_a^{a+\eps} h_1(x)\dd x =
  \frac{1}{\eps}\cdot\frac{\eps}{4}=\frac{1}{4}. 
\end{align}
However, by the Lebesgue differentiation theorem, as $\eps\to0$, this
converges \aex{} to $h_1(x)$. Hence, 
\begin{equation}
  \label{h1=}
h_1(x)=\frac{1}4 \qquad \text{a.e. $x\in\oi$}.
\end{equation}

We may now complete the proof.
It follows from \eqref{h1=} that 
$h_1(\psi(x))=\frac{1}{4}$ \aex{} on $\gO$, while
\eqref{h} implies that $h(x)\neq\frac{1}4$ \aex{} on $\oi$,
and thus
$h(\gf(x))\neq\frac14$ \aex{} on $\gO$.
Thus \eqref{hh} yields a contradiction.

Consequently, there is no graphon $W_1$ equivalent to $W$ with increasing
degree function.
\qed


\newcommand\AAP{\emph{Adv. Appl. Probab.} }
\newcommand\JAP{\emph{J. Appl. Probab.} }
\newcommand\JAMS{\emph{J. \AMS} }
\newcommand\MAMS{\emph{Memoirs \AMS} }
\newcommand\PAMS{\emph{Proc. \AMS} }
\newcommand\TAMS{\emph{Trans. \AMS} }
\newcommand\AnnMS{\emph{Ann. Math. Statist.} }
\newcommand\AnnPr{\emph{Ann. Probab.} }
\newcommand\CPC{\emph{Combin. Probab. Comput.} }
\newcommand\JMAA{\emph{J. Math. Anal. Appl.} }
\newcommand\RSA{\emph{Random Struct. Alg.} }
\newcommand\ZW{\emph{Z. Wahrsch. Verw. Gebiete} }
\newcommand\DMTCS{\jour{Discr. Math. Theor. Comput. Sci.} }

\newcommand\AMS{Amer. Math. Soc.}
\newcommand\Springer{Springer-Verlag}
\newcommand\Wiley{Wiley}

\newcommand\vol{\textbf}
\newcommand\jour{\emph}
\newcommand\book{\emph}
\newcommand\inbook{\emph}
\def\no#1#2,{\unskip#2, no. #1,} 
\newcommand\toappear{\unskip, to appear}

\newcommand\arxiv[1]{\texttt{arXiv:#1}}
\newcommand\arXiv{\arxiv}

\def\nobibitem#1\par{}


\begin{thebibliography}{99}

\bibitem{BR}  
B{\'e}la Bollob{\'a}s  \& Oliver Riordan:
Metrics for sparse graphs.
\emph{Surveys in Combinatorics 2009}, 211--287,
London Math. Soc. Lecture Note Ser. 365,
Cambridge Univ. Press, Cambridge, 2009. 

\bibitem{BCL}
Christian Borgs,
Jennifer Chayes \&
L\'aszl\'o Lov\'asz:
Moments of two-variable functions and the uniqueness of graph limits. 
\emph{Geom. Funct. Anal.} 19 (2010), no. 6, 1597--1619. 


  
\bibitem{SJ238} 
Persi Diaconis, Susan Holmes \& Svante Janson:
Threshold graph limits and random threshold graphs. 
\emph{Internet Mathematics} \vol5 (2008), no. 3, 267--318.  

\bibitem{SJ209}
Persi Diaconis \& Svante Janson:
Graph limits and exchangeable random graphs. 
\emph{Rend. Mat. Appl. (7)} \vol{28} (2008), no. 1, 33--61. 

\bibitem{SJ249}  
Svante Janson:
Graphons, cut norm and distance, couplings and rearrangements. 
\emph{NYJM Monographs} \vol4, 2013.  

\bibitem[Lov\'asz(2010)]{Lovasz}
L\'aszl\'o Lov\'asz:
\emph{Large Networks and Graph Limits},
American Mathematical Society, Providence, RI, 2012.

\bibitem{LovaszSz}
L\'aszl\'o Lov\'asz and Bal\'azs Szegedy:
Limits of dense graph sequences.
\emph{J. Combin. Theory Ser. B} \vol{96} (2006), no. 6, 933--957. 
  



\end{thebibliography}
\end{document}